\documentclass[preprint,12pt]{elsarticle}

\usepackage{amssymb}

\newtheorem{thm}{Theorem}[section]

\newtheorem{prp}[thm]{Proposition}

\newproof{pf}{Proof}
\newdefinition{df}[thm]{Definition}
\newdefinition{exm}[thm]{Example}

\newcommand{\NN}{{\mathbb{N}}}

\newcommand{\RR}{{\mathbb{R}}}

\newcommand{\cC}{{\cal C}}
\newcommand{\cD}{{\cal D}}
\newcommand{\tr}{{\rm tr \,}}
\newcommand{\re}{{\rm Re \,}}

\begin{document}

\begin{frontmatter}

\title{The spread of the spectrum of a nonnegative matrix with a zero diagonal element\tnoteref{label1}}
\tnotetext[label1]{The paper will appear in Linear Algebra and its Applications.}

\author{Roman Drnov\v{s}ek}
\ead{roman.drnovsek@fmf.uni-lj.si}


\address{Department of Mathematics, Faculty of Mathematics and Physics, University of Ljubljana, 
Jadranska 19, SI-1000 Ljubljana, Slovenia}

\baselineskip 7mm

\begin{abstract}
Let $A = [a_{i j}]_{i,j=1}^n$ be a nonnegative matrix with $a_{1 1} = 0$.
We prove some lower bounds for the spread $s(A)$ of $A$ that is defined as the maximum distance between any two eigenvalues of $A$.
If $A$ has only two distinct eigenvalues, then $s(A) \ge \frac{n}{2(n-1)} \, r(A)$, 
where $r(A)$ is the spectral radius of $A$. Moreover, this lower bound is the best possible.
\end{abstract}


\end{frontmatter}

\baselineskip 7mm



\section{Introduction}

Let $A$ be a complex $n \times n$ matrix with the spectrum $\{\lambda_1, \lambda_2, \ldots, \lambda_n\}$.
The spectral radius and the trace of $A$ are denoted by $r(A)$ and $\tr (A)$, respectively.
The spread $s(A)$ of $A$ is the maximum distance between any two eigenvalues, that is,
$s(A) = \max_{i,j} | \lambda_i - \lambda_j|$. This quantity was introduced by Mirsky \cite{Mir}, and it has been studied by several authors;
see e.g. \cite{MK} and the references therein. Note that $s(\lambda A)=  |\lambda| s(A)$ for every complex number $\lambda$ and 
that the spread of a nilpotent matrix is zero. Thus, when studying the spread of a matrix $A$, there is no loss of generality in assuming 
that $r(A) = 1$. 

Let $\cC_n$ (with $n \ge 2$) be the collection of all nonnegative $n \times n$ matrices $A = [a_{i j}]_{i,j=1}^n$ such that 
$a_{1 1} = 0$ and $r(A) = 1$.  It is not difficult to prove (see e.g. Proposition \ref{easy}) that the spread of a matrix $A \in \cC_n$ 
cannot be zero, that is, the number $1$ cannot be the only point in the spectrum of $A$. 
This motivates searching for lower bounds for the spread of $A$.  If $A$ has only two distinct eigenvalues, 
we prove that $s(A) \ge \frac{n}{2(n-1)}$, and we provide a matrix for which this lower bound is achieved.
Such a matrix is necessarily irreducible, that is, there exists no permutation matrix $P$ such that
$$ P^T A P = \left[\matrix{ 
A_{1 1} &  A_{1 2} \cr
   0    &  A_{2 2} }
\right] , $$
where $A_{1 1}$ and $A_{2 2}$ are square matrices.
\\

\section{Results}

We start with an easy observation.

\begin{prp}
\label{easy}
Let $A$ be a nonnegative $n \times n$ matrix with the spectral radius $r(A) = 1$. If $A$ has $k$ zero diagonal elements, then 
$$ s(A) \geq \frac{k}{n} \ . $$ 
In particular, if $A \in \cC_n$ then 
$$ s(A) \geq \frac{1}{n} \ . $$
\end{prp}

\begin{pf} 
Since $A$ is a nonnegative matrix, the spectral radius $r(A) = 1$ is its Perron eigenvalue. We denote it by $\lambda_1$,
while the rest eigenvalues of $A$ are denoted by $\lambda_2$, $\lambda_3$, $\ldots$, $\lambda_n$. 
For every $i=1,2, \ldots, n$ we have 
$$ \re (1- \lambda_i) \le |1- \lambda_i| = |\lambda_1- \lambda_i| \le s(A) , $$
and so $1 - s(A) \le \re \lambda_i$. It follows that 
$$ n (1 - s(A)) \le \sum_{i=1}^n \re \lambda_i = \sum_{i=1}^n \lambda_i = \tr (A) . $$
However, $\tr (A) = \sum_{i=1}^n a_{i i} \le n - k$, as $A$ has $k$ zero diagonal elements and 
$a_{i i} \le r(A) = 1$ for all $i$. We thus obtain that $n (1 - s(A)) \le  n - k$, and so $n \, s(A) \ge k$ as asserted.
\qed \end{pf}

Applying the known inequalities of Johnson, Loewy and London we will prove a better result for matrices in $\cC_n$.  
Let $A$ be a nonnegative $n \times n$ matrix and let $s_k := \tr (A^k)$ for $k \in \NN$. 
The JLL-inequalities (discovered independently by Loewy and London \cite{LL},  and Johnson \cite{Jo}) state that
$$ s_k^m \le n^{m-1} s_{k m} $$
for all positive integers $k$ and $m$. 
A slight modification of their proof gives the following inequalities.

\begin{prp}
\label{JLL}
Let $A$ be a nonnegative $n \times n$ matrix with $k$ zero diagonal elements. Then 
$$ s_1^m \le (n-k)^{m-1} s_{m} $$
for all $m \in \NN$. 
In particular, if $A \in \cC_n$ then 
$$ s_1^m \le (n-1)^{m-1} s_{m} $$
for all $m \in \NN$. 
\end{prp}

\begin{pf} 
Since $A$ is a nonnegative matrix, we have
$$ s_m = \tr (A^m) \ge \sum_{i=1}^n a_{i i}^m = \sum_{i \in J} a_{i i}^m , $$
where $J = \{i \in \{1, 2, \ldots, n\}: a_{i i} > 0\}$. 
On the other hand, H\"{o}lder's inequality gives 
$$ s_1^m = \left( \sum_{i \in J} a_{i i} \right)^m \le (n-k)^{m-1} \sum_{i \in J} a_{i i}^m , $$
and so we conclude that $s_1^m \le (n-k)^{m-1} s_{m}$.
\qed \end{pf}

Using Proposition \ref{JLL} we prove the following lower estimates for the spread of a matrix in $\cC_n$. 

\begin{thm}
If $A \in \cC_n$ then 
$$ s(A) > \frac{2}{4 + \sqrt{2 (n + 3)}} $$ 
for $n \ge 6$,  
$$ s(A) \ge \frac{5}{8 + \sqrt{74}} $$
for $n=5$, and 
$$ s(A) \ge \frac{1}{3} $$ 
for $n=4$. 
\end{thm}

\begin{pf} 
Since $s(A) > 0$ by Proposition \ref{easy} and since the result is true if $s(A) \ge 1$, 
we may assume that $s := s(A) \in (0, 1)$, and consequently the eigenvalues of $A$ have positive real parts.
Let $\lambda_1 = r(A)= 1$, $\lambda_2$, $\lambda_3$, $\ldots$, $\lambda_n$ be the spectrum of $A$.
By Proposition \ref{JLL}, we have 
$$ \left( \sum_{i=1}^n \lambda_i \right)^2 = s_1^2 \le (n-1) s_2 = (n-1) \sum_{i=1}^n \lambda_i^2 . $$
This inequality can be rewritten in the form
\begin{equation}
\sum_{i=1}^n \lambda_i^2 \le \sum_{i=1}^{n-1} \sum_{j=i+1}^n (\lambda_i - \lambda_j)^2 . 
\label{equiv}
\end{equation}
The right-hand side of (\ref{equiv}) is clearly at most $n(n-1) s^2/2$.
To obtain a lower bound for the left-hand side of (\ref{equiv}), we choose any eigenvalue $\lambda$ of $A$. 
Since $\lambda + \overline{\lambda} = 2 \, \re \lambda \ge 2(1-s) > 0$, we have  
$$ \lambda^2 + \overline{\lambda}^2 = (\lambda + \overline{\lambda})^2 - 2 |\lambda|^2 \ge (2(1-s))^2 - 2 = 4 s^2 - 8 s + 2 ,$$
and so we obtain the following lower bound for the left-hand side of (\ref{equiv}):
$$ \sum_{i=1}^n \lambda_i^2 = 1 + \sum_{i=2}^n \lambda_i^2 \ge 1 + \frac{n-1}{2}(4 s^2 - 8 s + 2) . $$
Therefore, the inequality (\ref{equiv}) gives the inequality
$$ \frac{n(n-1)}{2} s^2 \ge 1 + \frac{n-1}{2}(4 s^2 - 8 s + 2) , $$
which leads to the inequality
\begin{equation}
(n-1)(n-4) s^2 + 8(n-1)s - 2 n \ge 0 . 
\label{eq}
\end{equation}
For $n=4$ we obtain that $s \ge \frac{1}{3}$, while for $n=5$ we have 
$$ 2 s^2 + 16 s - 5 \ge 0 , $$
implying that 
$$ s \ge \frac{-8 + \sqrt{74}}{2} = \frac{5}{8 + \sqrt{74}} . $$ 
If $n \ge 6$ we rewrite the inequality (\ref{eq}) to the form 
$$ (n^2 - 5 n) s^2 + 8 n s - 2 n \ge -4 s^2 + 8 s = 4 s (2 - s) > 0 , $$ 
and so 
$$ (n - 5) s^2 + 8 s - 2 > 0 .  $$ 
It follows that 
$$ s > \frac{-4 + \sqrt{2 (n + 3)}}{n-5} = \frac{2}{4 + \sqrt{2 (n + 3)}} . $$ 
This completes the proof.
\qed \end{pf}

For $n \in \{2,3\}$ we can obtain sharp lower bounds for the spread of a matrix in $\cC_n$. 

\begin{prp}
If $A \in \cC_2$ then $s(A) \ge 1$; if $A \in \cC_3$ then $s(A) \ge \frac{3}{4}$.
Both bounds are exact.
\end{prp}

\begin{pf} 
Let $1$ and $\lambda$ be the eigenvalues of $A \in \cC_2$. By Proposition \ref{JLL}, we have 
$$ (1+\lambda)^2 = s_1^2 \le s_2 = 1 + \lambda^2 , $$
and so $\lambda  \le 0$ proving that $s(A) \ge 1$. The diagonal matrix ${\rm diag \,} (0,1) \in \cC_2$ shows that 
this lower bound is exact.

In the case $n=3$ we first suppose that a matrix $A \in \cC_3$ has real eigenvalues $1$, $\lambda$ and $\mu$.
We may assume that $0 \le \lambda \le \mu \le 1$. Then the inequality (\ref{equiv}) gives the inequality
$$ 1 + \lambda^2 + \mu^2 \le (1-\lambda)^2 + (1-\mu)^2 + (\lambda-\mu)^2 , $$
and so 
$$ 2 \lambda^2 \le 2 \lambda \mu \le (1-\lambda)^2 + (1-\mu)^2 - 1 \le 2 (1-\lambda)^2 - 1 = 2 \lambda^2 - 4 \lambda + 1 . $$
It follows that $\lambda \le \frac{1}{4}$, so that $s(A) \ge \frac{3}{4}$.

Assume now that a matrix $A \in \cC_3$ has eigenvalues $1$, $\lambda = a + i b$ and $\overline{\lambda} = a - i b$,
where $a \in \RR$ and $b > 0$. By Proposition \ref{JLL}, we have 
$$ (1 + 2 a)^2 = s_1^2 \le 2 s_2 = 2 (1 + \lambda^2 + \overline{\lambda}^2) = 2 + 4 a^2 - 4 b^2 \le 2 + 4 a^2, $$
and so $a \le \frac{1}{4}$. This implies that $s(A) \ge \frac{3}{4}$ as asserted.

The exactness of this lower bound is proved by the matrix 
$$ A = \frac{1}{4} 
\left[\matrix{ 
0 & 2 & 0 \cr
0 & 3 & 1 \cr
2 & 0 & 3 }
\right] \in \cC_3 $$
the spectrum of which is $\{1, \frac{1}{4}, \frac{1}{4}\}$.  
\qed \end{pf}

For $n \ge 4$ it looks difficult to obtain exact lower bounds for the spread of matrices in $\cC_n$. 
We thus restrict our attention to a special subset of $\cC_n$.
Proposition \ref{easy} trivially implies that every matrix in $\cC_n$ has at least two distinct eigenvalues, that is, 
$1$ is not the only point in its spectrum. Let $\cD_n$ (with $n \ge 2$) be the collection of all matrices in $\cC_n$ having exactly two distinct eigenvalues. We now prove sharp lower bounds for the spread of matrices in $\cD_n$.

\begin{thm}
If $A \in \cD_n$ then 
$$ s(A) \ge \frac{n}{2(n-1)} $$
Moreover, this bound is the best possible, i.e., there is a (necessarily irreducible) matrix $A \in \cD_n$
such that $s(A) = \frac{n}{2(n-1)}$.
\end{thm}

\begin{pf} 
Assume first that a matrix $A \in \cD_n$ is irreducible. Then $1$ is a simple eigenvalue of $A$ by the Perron-Frobenius theorem. 
Therefore, $A$ also has an eigenvalue $\lambda \in (-1, 1)$ of multiplicity $n-1$. 
In this case the inequality  (\ref{equiv}) reads as follows:
$$ 1 + (n-1) \lambda^2 \le (n-1) (1-\lambda)^2 .$$
Simplifying it, we obtain
$$ \lambda  \le \frac{n-2}{2(n-1)} . $$
This implies that 
$$ s(A) = 1 - \lambda \ge \frac{n}{2(n-1)} . $$

Assume now that a matrix $A \in \cD_n$ is reducible. 
Then, up to similarity with a permutation matrix, we may assume that
$$ A= \left[\matrix{ 
A_{11} & A_{12} & A_{13} & \ldots & A_{1m} \cr
   0   & A_{22} & A_{23} & \ldots & A_{2m} \cr
   0   &   0    & A_{33} & \ldots & A_{3m} \cr
\vdots & \vdots & \vdots & \ddots & \vdots \cr
   0   &   0    &   0    & \ldots & A_{mm} }
\right]  $$
where each of $A_{11}$, $A_{22}$, $\ldots$, $A_{mm}$ is either an irreducible (square) matrix or a $1 \times 1$ block.
Let $A_{k k}$ be one of these diagonal blocks that has a zero diagonal element. 
Without loss of generality we may assume that $s(A) < 1$, so that $0$ is not in the spectrum of $A$ implying that all $1 \times 1$ diagonal blocks are non-zero. Therefore, if $A_{k k}$ is an $r \times r$ matrix, then $r \ge 2$, and so
$$ s(A) \ge s(A_{k k}) \ge \frac{r}{2(r-1)} > \frac{n}{2(n-1)} . $$
This completes the proof of the first assertion of the theorem.

To show that the lower bound can be achieved, we define the matrix
$A = [a_{i,j}]_{i,j=1}^n$ with nonzero elements:
$a_{i,i+1} = n-i$ for $i=1, 2, \ldots, n-1$,
$a_{i,i} = n$  for $i=2, 3,  \ldots, n$, and 
$a_{i,j} = 2$ if $i - j$ is an even positive integer.
We also introduce the upper triangular matrix $U = [u_{i,j}]_{i,j=1}^n$ with nonzero elements:
$u_{i,i+1} = n-i$ for $i=1, 2, \ldots, n-1$, $u_{1,1} = 2(n-1)$ and 
$u_{i,i} = n-2$  for $i=2, 3,  \ldots, n$.
For example, if $n=5$ then 
$$ A = \left[\matrix{ 
0 & 4 & 0 & 0 & 0 \cr
0 & 5 & 3 & 0 & 0 \cr
2 & 0 & 5 & 2 & 0 \cr
0 & 2 & 0 & 5 & 1 \cr
2 & 0 & 2 & 0 & 5}
\right] 
\ \ \ \textrm{and} \ \ \
U = \left[\matrix{ 
8 & 4 & 0 & 0 & 0 \cr
0 & 3 & 3 & 0 & 0 \cr
0 & 0 & 3 & 2 & 0 \cr
0 & 0 & 0 & 3 & 1 \cr
0 & 0 & 0 & 0 & 3}
\right] .$$
The proof is complete if we show that $A$ and $U$ are similar matrices, because then we have 
$r(A) = 2(n-1)$, $s(A) = n$, and $\frac{1}{2(n-1)} A \in \cD_n$.
Define two nilpotent matrices
$$ N = \left[\matrix{ 
0 & 0 & 0 & 0 & \ldots & 0 & 0\cr
1 & 0 & 0 & 0 & \ldots & 0 & 0\cr
0 & 1 & 0 & 0 & \ldots & 0 & 0\cr
0 & 0 & 1 & 0 & \ldots & 0 & 0 \cr
\vdots & \vdots & \vdots & \vdots & \ddots & \vdots &\vdots \cr
0 & 0 & 0 & 0 & \ldots & 0 & 0 \cr 
0 & 0 & 0 & 0 & \ldots & 1 & 0}
\right] $$
and
$$ M = \left[\matrix{ 
0 & n\!-\!1 & 0 & 0 & \ldots & 0 & 0\cr
0 & 0 & n\!-\!2 & 0 & \ldots & 0 & 0\cr
0 & 0 & 0 & n\!-\!3 & \ldots & 0 & 0\cr
0 & 0 & 0 & 0 & \ldots & 0 & 0 \cr
\vdots & \vdots & \vdots & \vdots & \ddots & \vdots &\vdots \cr
0 & 0 & 0 & 0 & \ldots & 0 & 1 \cr 
0 & 0 & 0 & 0 & \ldots & 0 & 0}
\right] .$$
Introduce also the matrix 
$$ S = (I+N)(I-N)^{-1} = (I+N)(I+N+N^2+N^3+\ldots+N^{n-1}) = $$
$$ = I+2N+2N^2+2N^3+2N^4+\ldots+2N^{n-1} = $$
$$ = \left[\matrix{ 
1 & 0 & 0 & 0 & \ldots & 0 & 0\cr
2 & 1 & 0 & 0 & \ldots & 0 & 0\cr
2 & 2 & 1 & 0 & \ldots & 0 & 0\cr
2 & 2 & 2 & 1 & \ldots & 0 & 0 \cr
\vdots & \vdots & \vdots & \vdots & \ddots & \vdots &\vdots \cr
2 & 2 & 2 & 2 & \ldots & 1 & 0 \cr 
2 & 2 & 2 & 2 & \ldots & 2 & 1}
\right] .$$
Let $e_1$, $\ldots$, $e_n$ be the standard basis vectors, and let $e=e_1+\ldots+e_n=(1,1,\ldots,1)^T$. 
Observe that 
$$ A = M + nI - n e_1 e_1^T + 2 (N^2 + N^4 + N^6 + \ldots) = $$
$$ = M + (n-2)I - n e_1 e_1^T + 2 (I-N^2)^{-1} $$
and 
$$ U = M + (n-2)I + n e_1 e_1^T . $$
Note also that $[N,M] := N M - M N = I - n e_1 e_1^T$. By induction one can verify that 
$[N^k,M] = k N^{k-1} - n e_k e_1^T$ for $k=1,2, \ldots, n$. Then the commutator of $S$ and $M$ is 
$$ [S, M] = 2 \sum_{k=1}^n [N^k,M] = 2 \sum_{k=1}^n k N^{k-1} - 2 n \sum_{k=1}^n e_k e_1^T = 
2 (I-N)^{-2} - 2 n e e_1^T . $$
Now we have
$$ S U - A S = [S, M] +  n (S e_1) e_1^T + n e_1 e_1^T S -  2 (I-N^2)^{-1} S = $$
$$ = 2 (I-N)^{-2} - 2 n e e_1^T + n (2 e - e_1) e_1^T + n e_1 e_1^T - 2 (I-N^2)^{-1} (I+N)(I-N)^{-1} = $$
$$ = 2 (I-N)^{-2} - 2 (I-N)^{-2} = 0 . $$
This proves that the matrices $A$ and $U$ are similar. 
\qed \end{pf}

\vspace{5mm}
{\bf
\begin{center}
 Acknowledgments.
\end{center}
} 
The author was supported in part by the Slovenian Research Agency. 
He would like to thank Thomas Laffey and Helena \v{S}migoc for pointing out that Proposition \ref{JLL} holds. \\

\end{document}